\documentstyle[12pt]{article}
\begin{document}

{\noindent \small Dynamics of Continuous, Discrete and Impulsive
Systems, Series A: Mathematical Analysis, 11 (2004), 41-55.}

\begin{center}
{\large \bf PERIODIC SOLUTIONS FOR A CLASS OF\\ SINGULARLY
PERTURBED SYSTEMS $^\dagger$}
\end{center}
\begin{center}
Mikhail Kamenski$^1$, Oleg Makarenkov$^1$ and Paolo Nistri$^2$
\end{center}
\begin{center}
$^1$Department of Mathematics,\\
 Voronezh State University,
Voronezh, Russia.\\
$^2$ Dipartimento di Ingegneria dell'Informazione,\\
Universit\`a di Siena, 53100 Siena, Italy.
\end{center}

\date{}
\def\thefootnote{\fnsymbol{footnote}}
\footnotetext[2]{Supported by the research project ``Qualitative
analysis and control of dynamical systems'' at the University of Siena
and by RFRB grants 02-01-00189 and 02-01-00307.}
\def\thefootnote{\fnsymbol{footnote}}

\newtheorem{theorem}{Theorem}
\newtheorem{lemma}{Lemma}
\newtheorem{proposition}{Proposition}
\newtheorem{corollary}{Corollary}
\newtheorem{definition}{Definition}
\newtheorem{remark}{Remark}

\def\thedefinition{\arabic{definition}}
\def\thecorollary{\arabic{corollary}}
\def\thetheorem{\arabic{theorem}}
\def\vs{\vskip12pt}
\def\qed{\hfill $\squar$}
\def\squar{\vbox{\hrule\hbox{\vrule height 6pt \hskip
6pt\vrule}\hrule}}
\def\sqr{{\unskip\nobreak\hfil\penalty50\hskip2em\hbox{}\nobreak\hfil
{\squar}
\def\co{{\ol{\rm co}\,}}

        \parfillskip=0pt\finalhyphendemerits=0\par}}
\newenvironment{proof}{\noindent{\bf Proof}\ \ }{\sqr\par\vs}
\def\inte{\int\limits}     \let\qq=\qquad       \let\q=\quad
\let\w=\widetilde          \let\wh=\widehat     \let\mx=\mbox
\let\ol=\overline          \let\D=\Delta        \let\d=\delta
\let\e=\epsilon            \let\g=\gamma        \def\mm{{-\!\circ}}
\let\G=\Gamma              \let\ba=\beta        \def\dist{{\fam0 dist\,}}
\let\a=\alpha              \let\th=\theta       \let\nn=\nonumber
\let\s=\sigma              \let\O=\Omega      \def\L{{\cal L}}
\def\N{\mbox{\bf N}}       \def\bp{\mbox{\bf P}} \let\emp=\emptyset
\def\Z{\mbox{\bf Z}}       \def\ind#1{\mathop{#1}\limits}
\def\C{\mbox{\bf C}}       \let\l=\lambda       \let\sm=\setminus
\let\o=\omega              \def\meas{{\fam0 meas\,}}
\def\nor#1{{\left\|\,#1\,\right\|}}             \let\bcup=\bigcup
\def\R{\mbox{\bf R}}       \let\vf=\varphi    \def\sca#1#2{\langle #1,#2\rangle}
\def\Re{{\rm Re\,}}         \def\Im{{\rm Im\,}}        \def\DD{{\cal D}}
\let\ds=\displaystyle      \let\lra=\longrightarrow
\def\dfrac#1#2{\ds{#1\over #2}}\let\y=\eta
\let\sb=\subset            \def\J{{\cal J}}
\def\toto{{\begin{array}{c}\rightarrow\\*[-2ex]
             \rightarrow \end{array}}}
\catcode`\@=11
\def\section{\@startsection {section}{1}{\z@}{-3.5ex plus-1ex minus
-.2ex}{2.3ex plus.2ex}{\bf}}

\noindent {\bf Abstract.}$\;\;$ In this paper we provide
conditions to ensure the existence, for $\e>0$ sufficiently small,
of periodic solutions of given period $T>0$ in a prescribed domain
$U$ for a class of singularly perturbed first order differential
systems. Here $\e>0$ is the perturbation parameter. Our approach,
based on the topological degree theory and the averaging theory,
permits to weaken the conditions in (\cite{5}, Theorem~2) under
which the existence of periodic solutions is proved.

\vskip0.4truecm
\noindent
{\bf Keywords:} periodic solutions, averaging method, singular perturbations,
topological degree.

\vskip0.4truecm
\noindent
{\bf AMS (MOS) subject classification:}
34C25, 34C29, 34D15, 47H11.
\setcounter{section}{-1}

\vspace{6mm}
\noindent
{\bf Introduction}

\vspace{2mm}
\noindent
The starting point for the present work is the paper by K. Schneider
\cite{5}, devoted to the extension of the theory of vibrational stabilizability
to singularly perturbed first order control problems.
A basic tool for such extension is represented by the averaging theory.
In order to apply this theory the author assumes
that an appropriate coordinate trasformation of the fast variable is a
periodic diffeomorphism of fixed period. Indeed, this trasformation reduces
the considered system to the standard form for the application of the classical
averaging principle, (see, for instance, \cite{2}), this permits to prove the
existence
of periodic solutions for sufficiently small values of the perturbation
parameter
$\e>0$. The change of variable is introduced by means of a differential equation
which is supposed to have a $T$-periodic solution for any initial condition
in a suitable ball (assumption (A$_3$) of \cite{5}).

Following \cite{1}, \cite{4} and \cite{6} we propose here an
approach which combines the topological degree theory and the
averaging theory. This approach allows us to assume
the periodicity condition for the coordinate trasformation only on the
boundary of a given domain and then to derive, for any $\e>0$ sufficiently
small, the existence of a periodic solution for the system
whose fast component is contained in the interior of the domain.
We consider here the class of singularly perturbed systems
introduced by Schneider in his paper. However, in this paper we do
not look at the problem from the control point of view, but rather
we treat the general problem of the existence of periodic
solutions of singularly perturbed first order differential
systems. The main result together with the preliminaries and
assumptions are presented in Section 1. In Section 2 we illustrate
the main result by an example of singularly perturbed system
in $\R^3$ in which the boundary of
the domain is a circumference, which represents the periodic
solution of the first (fast) equation of the system at $\e=0$.

Finally, we would like to remark that as a direct consequence of our result we
obtain a result of the classical theory of ordinary differential equation
(\cite{3}, Theorem 3.1, p. 362).

\vspace{6mm}
\noindent
{\bf 1. Assumptions and Results}

\vspace{2mm}
\noindent
In this paper we consider the following system of differential equations
\begin{equation}\label{1}
\left\{\begin{array}{ll}
\dot x(t) = \e \phi(t, x(t), y(t))+\psi_1(t,x(t)),\\
\dot y(t) = \psi_2 (t,x(t),y(t)) - Ay(t),
\end{array}\right.
\end{equation}
where $\phi: \R\times \R^k\times \R^m \to \R^k$, $\psi_1: \R\times \R^k \to
\R^k$,
and $\psi_2:\R\times\R^k\times\R^m\to \R^m$ are continuous functions,
$T$-periodic with respect to time $t$. Moreover, $\psi_1$ is continuously
differentiable
with respect to the second variable $x$. $A$ is a
$m\times m$ matrix. We assume that the matrix $A$ has not
eigenvalues on the immaginary axis. Consequently the space $\R^m$
can be represented as the direct sum $E_+\oplus E_{-}$ of the
eigenspaces $E_{\pm}$ corresponding to the eigenvalues with
positive and negative parts respectively. If we denote by $A_+$
and $A_-$ the restrictions of $A$ to $E_+$ and $E_{-}$
respectively and by  $P_{+}$ and $P_{-}$ the projectors on $E_{+}$
and $E_{-}$ then we have the following dichotomy

\begin{equation}\label{2}
\|e^{-A_+t} P_+\| \le ce^{-\delta t}, \qquad \|e^{A_-t} P_-\| \le
ce^{-\delta t}, \qquad t\ge 0,
\end{equation}
for some $\delta>0$.
\noindent We introduce now in $\R^m$ the norm
given by
$$
\|x\|=\mbox{max}\{(\|P_+ x\|, \|P_- x\|)\}.
$$
We also assume that for every bounded set $B\subset \R^k$ there exists a
constant
$M(B)$ such that

\begin{equation}\label{3}
\|\psi_2(t,x,y)\|\le M(B) + \gamma \|y\|
\end{equation}
where $x\in B$ and

\begin{equation}\label{4}
\gamma<\frac{\delta}{c}.
\end{equation}
We denote by $\Omega(t, t_0, \xi)$ the solution of the Cauchy problem

\begin{equation}\label{5}
\left\{\begin{array}{ll}
\dot x(t) = \psi_1(t,x(t)), \qquad t\in [0,T],\\
x(t_0) = \xi.
\end{array}\right.
\end{equation}
We state the following property for the solution map $\xi\to\O(t, t_0, \xi)$,
which will be useful in the sequel.

\begin{lemma}
Let $z\in C^1([0,T], \R^m), \; f\in C([0,T], \R^m)$ and $b\in \R^m$. If

\begin{equation}\label{6}
\inte_0^t \frac{\partial \O}{\partial z}(s,0,z(s)) \dot z(s) ds + z(0)=
b + \inte_0^t f(s)ds,
\end{equation}
then
\begin{equation}\label{7}
z(t)= b + \inte_0^t \frac{\partial \O}{\partial z}(0,s, \O(s,0, z(s)))f(s) ds.
\end{equation}
\end{lemma}

\noindent
{\it Proof.} Take the derivative of (\ref{6}) with respect to $t$ obtaining
\begin{equation}\label{8}
\frac{\partial \O}{\partial z}(t,0,z(t)) \dot z(t) = f(t).
\end{equation}
Derive with respect to $\xi$ the identity
$$
\O(0,t, \O(t,0, \xi))=\xi
$$
to obtain
\begin{equation}\label{9}
\frac{\partial \O}{\partial z}(0,t,\O(t,0,\xi)) \frac{\partial \O}{\partial z}
(t,0,\xi) = I,
\end{equation}
whenever $\xi\in\R^k$. Applying now
$\dfrac{\partial \O}{\partial z}(0,t,\O(t,0,z(t)))$
to (\ref{8}) we get
$$
\dot z(t)= \frac{\partial \O}{\partial z}(0,t,\O(t,0,z(t))) f(t).
$$
Finally, integrating this equation from $0$ to $t$ and observing that from
(\ref{6}) we have $z(0)=b$, one has (\ref{7}).
\qed

\vskip0.5truecm
Let us now denote by $\eta_y(t,s,\xi)$ the solution of the Cauchy problem

\begin{equation}\label{10}
\left\{\begin{array}{ll}
\dot z(t) = \dfrac{\partial \psi_1}{\partial x}(t,\O(t,0,\xi))z(t)
+\phi(t,\O(t,0,\xi),y(t)),\\
z(s) = 0.
\end{array}\right.
\end{equation}
where $y$ is a $T$-periodic continuous function and $t,s\in[0,T]$. For
$y=0$ we write $\eta_0(t,s,\xi)$. In the sequel we will denote by
$\|y\|_{C_T}$ the norm of $y$ in the Banach space $C_T([0,T], \R^m)$
of $T$-periodic continuous functions.

\vskip0.5truecm
We can state the following main result.
\begin{theorem} Assume that there exists an open bounded set $U\subset\R^k$
such that

\begin{enumerate}
\item[(A$_1$)] $\O(T,0,\xi)=\xi$ for all $\xi\in \partial U$;

\item[(A$_2$)] $\eta_y(0,s,\xi)\not= \eta_y(T,s,\xi)$ for all $s\in[0,T]$,
$\xi\in\partial U$ and  $\|y\|_{C_T}\le \dfrac{cM}{\delta-\gamma c}$, where
$M=\mbox{sup}_{t\in[0,T]}\mbox{sup}_{\xi\in\overline U} \|\O(t,0,\xi)\|$;

\item[(A$_3$)] $\mbox{deg}\,(\eta_0(T,T,\cdot)-\eta_0(0,T,\cdot),U,0)\not= 0$.
\end{enumerate}
Then there exists $\e_0>0$ such that for $\e\in(0,\e_0)$ system (\ref{1})
has at least one $T$-periodic solution $(x,y)$ such that
$$
\O(0,t,x(t))\in U \quad \mbox{for any}\quad t\in[0,T] \qquad \mbox{and} \qquad
\|y\|_{C_T}\le
\frac{cM}{\delta-\gamma c}.
$$

\end{theorem}

\vspace{2mm} For the proof of this theorem we need the following
lemma.

\begin{lemma} Let $\eta(t,s)$ be the solution of the problem
$$
\left\{\begin{array}{ll}
\dot z(t) = A(t)z(t)+f(t),\\
z(s) = 0.
\end{array}\right.
$$
where $A(t)$ is a $T$-periodic continuous $k\times k$ matrix and $f$ is a
$T$-periodic
continuous function taking values in $\R^k$. We have that
$$
\eta(T,s)-\eta(0,s)=\inte_{s-T}^s X^{-1}(\tau) f(\tau)d\tau
$$
where $X(t)$ is a fundamental matrix of the linear system
$$
\dot z(t) = A(t)z(t),
$$
such that $X(0)=I.$
\end{lemma}

\noindent
{\it Proof of Lemma 2.} Observe that $\eta(t,s)=\inte_s^t
X(t)X^{-1}(\tau)f(\tau)d\tau$.
Therefore
\begin{equation}\label{11}
\eta(T,s)-\eta(0,s)=\inte_s^T X(T)X^{-1}(\tau) f(\tau)d\tau - \inte_s^0
X^{-1}(\tau) f(\tau)d\tau.
\end{equation}
Make the change of variable $\tau=u+T$ in the integral
$$
J=\inte_s^T X(T)X^{-1}(\tau) f(\tau)d\tau
$$
obtaining
$$
J=\inte_{s-T}^0 X(T)X^{-1}(T+u) f(T+u)du =
\inte_{s-T}^0 \Phi(T)e^{\Lambda T}X^{-1}(T+u) f(u)du=
$$
$$
=\inte_{s-T}^0 \Phi(T)e^{-\Lambda u}e^{\Lambda (T+u)}X^{-1}(T+u) f(u)du
$$
where $\Phi(t)e^{\Lambda t}$ is a Floquet representation of $X(t)$.
Therefore $e^{\Lambda t}X^{-1}(t)=\Phi^{-1}(t),$ then
$$
e^{\Lambda (T+u)}X^{-1}(T+u)= e^{\Lambda u}X^{-1}(u)
$$
and
$$
J=\inte_{s-T}^0 \Phi(0)e^{-\Lambda u}e^{\Lambda u}X^{-1}(u) f(u)du=
\inte_{s-T}^0 X^{-1}(u) f(u)du.
$$
Therefore from (\ref{11}) we have
\begin{equation}\label{12}
\eta(T,s)-\eta(0,s)=\inte_{s-T}^0 X^{-1}(\tau) f(\tau)d\tau + \inte_0^s
X^{-1}(\tau)f(\tau)d\tau=
$$
$$
=\inte_{s-T}^s X^{-1}(\tau) f(\tau)d\tau.
\end{equation}
Which is the claim.\qed

\begin{remark} From (\ref{9}), (\ref{12}) and the fact that $\dfrac{\partial
\O}{\partial z}(t,0,\xi)$
is the fundamental matrix of (\ref{10}) with $\phi=0$, we can deduce the
following formula
$$
\eta_0(T,T,\xi)-\eta_0(0,T,\xi)=\inte_0^T \frac{\partial \O}{\partial z}
(0,\tau,\O(\tau,0,\xi))\phi(\tau,\O(\tau,0,\xi),0)d\tau
$$
for $\xi\in \partial U$.
\end{remark}

\noindent
{\it Proof of Theorem 1.} The existence of $T$-periodic solutions for
system (\ref{1}) is equivalent to the existence of solution pairs
$(x,y)\in C_T([0,T], \R^k)\times C_T([0,T], \R^m)$ of the following system
\begin{equation}\label{13}
\left\{\begin{array}{lll}
x(t)=x(T) + \e \inte_0^t \phi(s,x(s),y(s))ds + \inte_0^t \psi_1(s,x(s))ds, \\
y_+ (t)=e^{ - A_+ t}(I-e^{- A_+ T})^{-1}\inte_0^T e^{- A_+ (T-s)}
P_+ \psi_2(s,x(s),y(s))ds +\\\quad\quad\quad\; + \inte_0^t e^{-
A_+ (t-s)}
P_+ \psi_2(s,x(s),y(s))ds ,\\
y_- (t)=e^{ A_- (T-t)}(e^{ A_- T}-I)^{-1}\inte_0^T e^{A_-s} P_-
\psi_2(s,x(s),y(s))ds -\\\quad\quad\quad\; - \inte_t^T e^{A_-
(s-t)} P_- \psi_2(s,x(s),y(s))ds,
\end{array}\right.
\end{equation}

\noindent
where $y_+=P_+ y$ and $y_-=P_- y.$ Consider the change of variable
\begin{equation}\label{14}
x(t)=\O(t,0,z(t)), \qquad t\in[0,T],
\end{equation}
with inverse given by
\begin{equation}\label{15}
z(t)=\O(0,t,x(t)), \qquad t\in[0,T].
\end{equation}
Observe that if $x$ is the first coordinate of the solution of system
(\ref{13}), then $x$ is differentiable and therefore from (\ref{14})
$z$ is also differentiable. Consider
\begin{equation}\label{16}
\dfrac{d}{dt} \O(t,0,z(t)) = \dfrac{\partial \O}{\partial t}(t,0,z(t))+
\dfrac{\partial \O}{\partial z}(t,0,z(t))\dot z(t),
\end{equation}
since
$$
\dfrac{\partial \O}{\partial t}(t,0,z(t))=\psi_1(t,\O(t,0,z(t)))
$$
from (\ref{16}) we have that
$$
\O(t,0,z(t)) - z(0) = \inte_0^t \psi_1(s,\O(s,0,z(s)))ds +
\inte_0^t \dfrac{\partial \O}{\partial z}(s,0,z(s))\dot z(s)ds,
$$
or equivalently,
\begin{equation}\label{17}
\O(t,0,z(t)) -  \inte_0^t \psi_1(s,\O(s,0,z(s)))ds =
z(0) +  \inte_0^t \dfrac{\partial \O}{\partial z}(s,0,z(s))\dot z(s)ds.
\end{equation}
Let $(x,y)$ be a solution of system (\ref{13}); by using (\ref{14}),
(\ref{17}) and (\ref{7}) with $b=\O(T,0,z(T))$ we can rewrite (\ref{13}) in the
following form
\begin{equation}\label{18}
\left\{\begin{array}{lll}
z(t)= \O(T,0,z(T)) + \e \inte_0^t \Phi(s,x(s),y(s))ds, \\
y_+ (t)=e^{ - A_+ t}(I-e^{- A_+ T})^{-1}
\cdot\\\quad\quad\quad\;\cdot
 \inte_0^T e^{- A_+ (T-s)}P_+ \psi_2(s,\O(s,0,z(s)),y(s))ds +
\\\quad\quad\quad\;
+\inte_0^t e^{- A_+ (t-s)}
P_+ \psi_2(s,\O(s,0,z(s)),y(s))ds ,\\
y_- (t)=e^{ A_- (T-t)}(e^{ A_- T}-I)^{-1}
\cdot\\\quad\quad\quad\;\cdot \inte_0^T e^{A_- s}P_-
\psi_2(s,\O(s,0,z(s)),y(s))ds - \\\quad\quad\quad\; - \inte_t^T
e^{A_- (s-t)}P_- \psi_2(s,\O(s,0,z(s)),y(s))ds,
\end{array}\right.
\end{equation}
where $\Phi(\tau,\xi,y)=\dfrac{\partial\O}{\partial
z}(0,\tau,\O(\tau,0,\xi))\phi(\tau,\O(\tau,0,\xi),y)$.
Therefore the problem of finding $T$-periodic solutions for (\ref{1}) is
equivalent to the problem
of the existence of zeros for the compact vector field
$$
G_\e\left(\begin{array}{lll} z\\y_+\\y_-\end{array}\right)(t)=
\left(\begin{array}{lll} z(t)-F_1(\e,z,y_+
+y_-)(t)\\y_+(t)-F_2(z,y_+ + y_-)(t)
\\y_-(t)-F_3(z,y_+ + y_-)(t)\end{array}\right)
$$
with $F_1,F_2$ and $F_3$ defined as the righthand sides of the
three equations in (\ref{18}).

\vspace{2mm}
Consider now in $C_T([0,T], \R^k)\times C_T([0,T],
E_+)\times C_T([0,T], E_-)$ the open set
$$
V_U = Z_U\times B(0,r)\times B(0,r)
$$
where $Z_U = \{z\in C_T([0,T], \R^k):\;\, z(t)\in U,\; t\in [0,T]\}$ and
$r>\dfrac{cM}{\delta-\gamma c}$.
Consider the homotopy
$$
H_\e\left(\lambda, \left(\begin{array}{lll} z\\y_+\\y_-\end{array}
\right)\right)(t) =
$$
$$
=\left(\begin{array}{lll} z(t)- \O(T,0,z(T)) -\e
\inte_0^{\mu(\lambda, t)} \Phi(\tau, z(\tau), y_+(\tau) +
y_-(\tau))d\tau\\
\qquad \qquad \qquad y_+(t)-\lambda F_2(z,y_+ + y_-)(t)\\
\qquad \qquad \qquad y_-(t)-\lambda F_3(z,y_+ +
y_-)(t)\end{array}\right)
$$
with $\mu(\lambda, t)=\lambda t+(1-\lambda)T$ and $\lambda\in
[0,1]$. We prove in the sequel that, for sufficiently small
$\e>0$, $H_\e$ is an admissible homotopy on $\partial V_U$. We
argue by contradiction, hence we assume the existence of sequences
$$
\lambda_n \to \lambda_0, \quad \lambda_n\in[0,1], \quad \e_n\to 0, \quad \e_n >
0 \quad \mbox{and} \qquad \left(\begin{array}{lll} z_n\\
y_+^n\\y_-^n\end{array}\right)\in \partial V_U
$$
such that
\begin{equation}\label{19}
\left\{\begin{array}{lll} z_n(t) = \O(T,0,z_n(T)) + \e_n
\inte_0^{\mu(\lambda_n,\, t)} \Phi(\tau, z_n(\tau), y_+^n(\tau) +
y_-^n(\tau))d\tau\\y_+^n(t)=
\lambda_n F_2(z_n,y_+^n+y_-^n)(t) \\
y_-^n(t)=\lambda_n F_3(z_n,y_+^n+y_-^n)(t).
\end{array}
\right.
\end{equation}
We first show that
\begin{equation}\label{20}
\|y_+^n\|_{C_T}, \; \|y_-^n\|_{C_T}\le \dfrac{cM}{\delta-c\gamma}.
\end{equation}
For this consider
$$
y_+^n(t)=\lambda_n\inte_0^T \sum_{m=0}^\infty
e^{-A_+((m+1)T+t-s)}P_+ \psi_2 (s,\O(s,0,z_n(s)),y_n(s))ds+
$$
$$
+\lambda_n \inte_0^t e^{-A_+(t-s)}P_
+\psi_2(s,\O(s,0,z_n(s)),y_n(s))ds = \lambda_n\inte_ {-\infty} ^t
e^{-A_+ (t-s)}P_+\tilde{\psi}_2(s)ds;
$$
and
$$
y_-^n (t)= -\lambda_n \inte_0^T \sum_{m=0}^\infty
e^{A_-((m+1)T+s-t)} P_-\psi_2(s,\O(s,0,z_n(s)),y_n(s))ds-
$$
$$
-\lambda_n \inte_t^T e^{A_-(s-t)}P_-\psi_2(s,\O(s,0,z_n(s)),y_n(s))ds
= - \lambda_n \inte_t ^{+\infty} e^{A_-
(s-t)}P_-\tilde{\psi_2}(s)ds,
$$
where $y_n=y_+^n+y_-^n$ and $\tilde\psi_2$ denotes the
$T$-periodic extension from $[0,T]$ to $\R$ of the function
$\psi_2(s,\O(s,0,z_n(s)),y_n(s))$. By using (\ref{2}) and (\ref{3}) we
have that
$$
\|y_+^n(t)\| \le \dfrac{cM}{\delta} + \dfrac{c \gamma
}{\delta}\|y\|_{C_T};
$$
$$
\|y_-^n(t)\| \le \dfrac{cM}{\delta} + \dfrac{c \gamma
}{\delta}\|y\|_{C_T}.
$$
and so (\ref{20}). By our choice of $r$ we have that
$$
\left( \begin{array}{lll} y_+^n\\y_-^n\end{array}\right)\notin
\partial(B(0,r)\times B(0,r)).
$$
Therefore, it must be $z_n\in \partial Z_U$ and so, for any $n\in\N$, there
exists
$t_n\in[0,T]$ such that $z_n(t_n)\in\partial U$. Assumption (A$_1$) implies that
\begin{equation}\label{21}
z_n(t_n)\in \O(T,0,z_n(t_n)).
\end{equation}
Putting $t=T$ and then $t=t_n$ in the first equation of (\ref{19}) and
subtracting
the obtained equations one has

\begin{equation}\label{22}
z_n(T)-z_n(t_n)=\e_n \inte_{\mu(\lambda_n,\, t_n)}^T \Phi (\tau,
z_n(\tau), y_n(\tau))d\tau.
\end{equation}
Furthermore subtracting (\ref{21}) from the first equation of
(\ref{19}), where we have replaced $t$ by $T$ we obtain
$$
z_n(T)-z_n(t_n)= \O(T,0,z_n(T)) - \O(T,0,z_n(t_n)) +
\e_n \inte_0^T \Phi (\tau, z_n(\tau), y_n(\tau))d\tau.
$$
Since $\O$ is differentiable we can rewrite the last equality as follows
\begin{equation}\label{23}
\begin{array}{lll}\left(I-\dfrac{\partial\O}{\partial z}(T,0,z_n(t_n))\right)
(z_n(T)-z_n(t_n))=\\
=\e_n \inte_0^T \Phi (\tau, z_n(\tau), y_n(\tau))d\tau + o(z_n(t_n),
z_n(T)-z_n(t_n)),
\end{array}
\end{equation}
where $o(\xi, h)$ is such that
$$
\frac{o(\xi,h)}{\|h\|} \to 0 \qquad \mbox{as} \qquad \|h\| \to 0
$$
uniformly with respect to $\xi$ belonging to compact sets. Replacing (\ref{22})
into
(\ref{23}) and dividing by $\e_n > 0$ we obtain
\begin{equation}\label{24}
\begin{array}{lll}\left(I-\dfrac{\partial\O}{\partial z}(T,0,z_n(t_n))\right)
\inte_{\mu(\lambda_n,\, t_n)}^T
\Phi (\tau, z_n(\tau), y_n(\tau))d\tau =\\\\
=\inte_0^T \Phi (\tau, z_n(\tau), y_n(\tau))d\tau +
\dfrac{o(z_n(t_n), z_n(T)-z_n(t_n))} {\e_n}.
\end{array}
\end{equation}
From (\ref{22}) there exists a constant $C>0$ such that
$$
\|z_n(T)-z_n(t_n)\|\le C\e_n.
$$
Therefore
\begin{equation}\label{25}
\dfrac{o(z_n(t_n), z_n(T)-z_n(t_n))}
{\e_n} \to 0 \qquad \mbox{as} \qquad n\to \infty.
\end{equation}
On the other hand the operators $F_2$ and $F_3$ are compact, and
so we can, without loss of generality, assume that the sequences
$\{y_+^n\}$, $\{y_-^n\}$ and consequently $y_n=y_+^n + y_-^n$
converge. Let $y_n \to y_0$, thus
$$
\|y_0\|_{C_T}\le \dfrac{cM}{\delta-c\gamma}.
$$
Furthermore, we can also assume that $z_n(t_n) \to \xi_0$, hence  $\xi_0 \in
\partial U$.
Since
$$
\|z_n(t) - z_n(t_n)\|= \e_n \|\inte_{l_n} \Phi (\tau, z_n(\tau), y_n(\tau))d\tau
\|,
$$
where $l_n$ is the segment joining $t$ with $\mu(\lambda_n, t_n)$,
we obtain
$$
z_n(\tau) \to \xi_0 \qquad \mbox{as} \qquad n \to \infty
$$
uniformly with respect to $t$.

We are now in a position to pass to the limit in (\ref{24}) obtaining
\begin{equation}\label{26}
\left(I-\dfrac{\partial\O}{\partial z}(T,0,\xi_0))\right)
\inte_{t_0}^T \Phi (\tau, \xi_0, y_0(\tau))d\tau
= \inte_0^T \Phi (\tau, \xi_0, y_0(\tau))d\tau,
\end{equation}
where $t_0=\lim_{n\to \infty} \mu(\lambda_n, t_n)$, $t_0\in
[0,T]$. \vspace{2mm}
Using the following property for the
traslation operator $\dfrac{\partial\O}{\partial z}$:
\begin{equation}\label{27}
\dfrac{\partial\O}{\partial z}(T,0,\xi)
\dfrac{\partial\O}{\partial z}(0,\tau,\O(\tau,0,\xi))=
\dfrac{\partial\O}{\partial z}(T,\tau,\O(\tau,0,\xi)),
\end{equation}
we can rewrite (\ref{26}) as follows
$$
\inte_{t_0}^T \dfrac{\partial\O}{\partial z}(T,\tau,\O(\tau,0,\xi_0))
\phi(\tau,\O(\tau,0, \xi_0), y_0(\tau))d\tau + \inte_0^{t_0}
\Phi (\tau, \xi_0, y_0(\tau))d\tau=0,
$$
or equivalently, by the definition of $\eta_{y_0},$ we can write
\begin{equation}\label{28}
\eta_{y_0}(T,t_0,\xi_0) + \inte_0^{t_0}
\Phi (\tau, \xi_0, y_0(\tau))d\tau=0.
\end{equation}
But,
$$
\eta_{y_0}(0,t_0,\xi_0) = \inte_{t_0}^0
\Phi (\tau, \xi_0, y_0(\tau))d\tau,
$$
hence (\ref{28}) takes the form
$$
\eta_{y_0}(T,t_0,\xi_0) - \eta_{y_0}(0,t_0,\xi_0)=0,
$$
which is a contradiction with assumption (A$_2$).

\vspace{2mm}
In conclusion, we have proved that for all
$\e\in(0,\e_0)$ the function $H_\e$ is an admissible homotopy on
$\partial V_U$. For $\lambda=0$ we have
$$
H_\e\left(0, \left(\begin{array}{lll}
z\\y_+\\y_-\end{array}\right)\right)(t)= \left(\begin{array}{lll}
z(t)- \O(T,0,z(T)) -\e \inte_0^T \Phi(\tau,
z(\tau), y(\tau)))d\tau
\\ \qquad\qquad\qquad\qquad y_+(t)
\\  \qquad\qquad\qquad\qquad y_-(t)
\end{array}\right).
$$
By the reduction property of the topological degree, for any $\e\in(0,\e_0)$,
one has

\begin{equation}\label{29}
\begin{array}{lll}
\mbox{deg}(G_\e, V_U, 0) = \mbox{deg}(H_\e(0, \cdot), V_U, 0) =\\\\
\mbox{deg}(H_\e(0, \cdot), V_U\bigcap C_{\mbox{const}}([0,T],
\R^k)\times\{0\}\times\{0\}, 0)
= \mbox{deg}(\Phi_{0,\e}, U, 0),
\end{array}
\end{equation}
where $C_{\mbox{const}}([0,T], \R^k)$ denotes the space of constant functions
defined
on $[0,T]$ with values in $\R^k$ and
$$
\Phi_{0,\e}(\xi) = \xi - \O(T, 0, \xi) -\e\inte_0^T \Phi(\tau, \xi, 0)d\tau.
$$
By (A$_1$) we have
$$
\O(T, 0, \xi)=\xi \qquad \mbox{for} \qquad \xi\in\partial U.
$$
Therefore, for $\xi\in\partial U$, we get
$$
\Phi_{0,\e}(\xi) = -\e \inte_0^T \Phi(\tau, \xi, 0)d\tau.
$$
But
\begin{equation}\label{30}
\mbox{deg}(\Phi_{0,\e}, U, 0) = \mbox{deg}(\Phi_{0,1}, U, 0),
\end{equation}
and by Remark 1 we have
$$
\Phi_{0,1}(\xi) = \eta_0(T,T,\xi) - \eta_0(0,T,\xi)
$$
Now using (\ref{29}), (\ref{30}), assumption (A$_3$) and the homotopy invariance
of the topological degree we obtain
$$
\mbox{deg}(G_\e, V_U, 0) \not= 0.
$$
This concludes the proof of the theorem. \qed

\vspace{6mm}
\noindent
{\bf 2. An Example}

\vspace{2mm}
\noindent
In this Section we provide an example to illustrate our main result: Theorem~1.
\begin{equation}\label{31}
\left\{\begin{array}{ll}
\dot x(t) = \e \phi(t, x(t), y(t))+\psi_1(t,x(t)),\\
\dot y(t) = \|x(t)\| - a y(t), \quad a>0, \quad t\in[0,2\pi]
\end{array}\right.
\end{equation}
where $\psi_1(t,x)=g(t,x)\sigma(x)+\psi(x)$, with
$$
\psi(x) = \left(\begin{array}{ll} x_2+x_1(x_1^2+x_2^2 - 1)\\
-x_1+x_2(x_1^2+x_2^2-1) \end{array}\right)
$$
and $(x_1,x_2)=x$. Let $x_0=\left
(\begin{array}{ll}\sin\theta\\\cos\theta\end{array}\right)$. We assume that
$$
\phi:\R\times\R^2\times\R\to \R^2,\quad
g:\R\times\R^2\to\R^2,\quad \sigma:\R^2\to\R
$$
are continuous differentiable functions. Moreover, the functions $\phi$ and
$g$ are $2\pi$-periodic with respect to time $t$ and $\sigma$ and $g$ satisfy
the following conditions
\begin{equation}\label{32}
\sigma(x_0(\theta))=0 \qquad \mbox{for any} \qquad \theta\in[0,2\pi],
\end{equation}
\begin{equation}\label{33}
\dot\sigma(x_0(\theta))=\dfrac{\partial g}{\partial x}(t, x_0(\theta))=0
\qquad \mbox{for any} \qquad t,\theta\in[0,2\pi].
\end{equation}
Denote by $U\subset \R^2$ the interior of the unitary circle centered at the
origin. At $\e=0$ the first equation of system (\ref{31}) has the form
\begin{equation}\label{34}
\dot x = \psi_1(t,x)= g(t,x)\sigma(x)+\psi(x),
\end{equation}
We denote by $\O(t,t_0,\xi)$ the solution of (\ref{34}) satisfying
the initial condition $x(t_0)=\xi$. Observe that due to condition
(\ref{32}) $x_0$ is the solution of equation (\ref{34}) satisfying
the initial condition $x(0)=\xi$, whenever $\|\xi\|=1$. Therefore
condition (A$_1$) of Theorem 1 holds true with $T=2\pi$. To verify conditions
(A$_2$), (A$_3$) we calculate in the sequel the translation
operator $\dfrac{\partial\O}{\partial z} (t,\tau,\O(\tau,0,\xi))$,
where $\xi=\left
(\begin{array}{ll}\sin\theta\\\cos\theta\end{array}\right)$ and
$t,\tau,\theta\in[0,2\pi]$. For this observe that if
$\xi=\left(\begin{array}{ll} \sin\theta\\ \cos\theta
\end{array}\right)$ then by \cite{3} and (\ref{33}) we get
\begin{equation}\label{35}
\dfrac{\partial\O}{\partial z}(t,\tau,\O(\tau,0,\xi)) = Y(t,\tau,\theta),
\end{equation}
where $Y:\R\times\R\times\R \to L(\R^2, \R^2)$, the vector
space of the linear operators from $\R^2$ to $\R^2$, is the solution of the
following
problem
\begin{equation}\label{36}
\left\{\begin{array}{ll}
\dfrac{dY}{dt}(t,\tau,\theta)= \dfrac{d\psi}{dx}(\O(t,0,\xi))
Y(t,\tau,\theta),\\\\
Y(\tau,\tau,\theta)=I.
\end{array}\right.
\end{equation}
Therefore the problem is reduced to find $Y(t,\tau,\theta)$. To this aim
we show that the function
\begin{equation}\label{37}
K(t,\theta)= \left(\begin{array}{ll} \quad\cos(t+\theta) \quad e^{2t}
\sin(t+\theta)\\
-\sin(t+\theta) \quad e^{2t} \cos(t+\theta)\end{array}\right)
\end{equation}
satisfies the equation
\begin{equation}\label{38}
\dfrac{dK}{dt}(t,\theta)= \dfrac{d\psi}{dx}(\O(t,0,\xi)) K(t,\theta).
\end{equation}
In fact,
$$
\dfrac{d\psi}{dx}(\O(t,0,\xi)) = \dfrac{d\psi}{dx}(x_0(t+\theta))=
$$
$$
= \left(\begin{array}{ll} \quad\qquad2\sin^2(t+\theta) \quad\quad\quad \quad
2\sin(t+\theta)\cos(t+\theta)+1\\
2\sin(t+\theta)\cos(t+\theta)-1\quad\quad\quad\quad
2\cos^2(t+\theta)\end{array}\right)
$$
and
$$
\dfrac{dK}{dt}(t,\theta) =
\left(\begin{array}{ll} -\sin(t+\theta) \quad 2e^{2t}\sin(t+\theta)+
e^{2t}\cos(t+\theta)\\
-\cos(t+\theta) \quad 2e^{2t}\cos(t+\theta)-e^{2t}\sin(t+\theta)
\end{array}\right).$$
Thus the function
\begin{equation}\label{39}
Y(t,\tau,\theta)=K(t,\theta)(K(\tau ,\theta))^{-1}
\end{equation}
is the solution of system (\ref{36}).
Let $\xi=\left(\begin{array}{ll} \sin\theta\\ \cos\theta
\end{array}\right)$ and denote by $\eta_y(t,s,\xi)$, $t,s\in [0,2\pi]$,
the solution of
$$
\left\{\begin{array}{ll}
\dot z(t) = \dfrac{\partial
\psi_1}{\partial x}(t,\O(t,0,\xi))z(t) +\phi(t,\O(t,0,\xi),y(t)),\\\\
z(s) = 0.
\end{array}\right.
$$
corresponding to the $2\pi$-periodic continuous function $y$.
By Lemma 1, (\ref{35}), (\ref{38}) and (\ref{39})
we have that
$$
\langle \eta_y(2\pi, s, \xi) - \eta_y(0, s, \xi), \, \psi(\xi)\rangle=
$$
$$
= \langle \inte_{s-2\pi}^s \dfrac{\partial\O}{\partial z}(0,\tau,\O(\tau,0,\xi))
\phi(\tau, \O(\tau,0,\xi),y(\tau))d\tau, \, \psi(\xi)\rangle=
$$
$$
= \langle \inte_{s-2\pi}^s K(0,\theta)(K(\tau,\theta))^{-1}
\phi(\tau, x_0(t+\theta),y(\tau))d\tau, \, \dot x_0(\theta)\rangle=
$$
$$
=\inte_0^{2\pi}\langle \phi(\tau, x_0(\tau+\theta), y(\tau)), \, \dot
x_0(\tau+\theta)\rangle d\tau=
$$
$$
=\inte_0^{2\pi}\langle \phi(r-\theta, x_0(r), y(r-\theta)), \, \dot
x_0(r)\rangle dr, \qquad \theta\in [0,2\pi].
$$
Hence if
\begin{equation}\label{40}
\inte_0^{2\pi}\langle \phi(r-\theta, x_0(r), y(r-\theta)), \, \dot x_0(r)\rangle
dr
>0, \quad \theta\in [0,2\pi], \quad \|y\|_{C_T}\le \dfrac{1}{a},
\end{equation}
then we have
$$
\eta_y(2\pi, s, \xi) - \eta_y(0, s, \xi)\not= 0, \quad s\in [0,2\pi], \quad
\xi\in\partial U, \quad \|y\|_{C}\le \dfrac{1}{a},
$$
and so assumption (A$_2$) is verified.
Moreover, using the property $\mbox{ind}(\psi, U, 0)=1$ we obtain
$$
\mbox{ind}(\eta_0(2\pi, 2\pi, \cdot) - \eta_y(0, 2\pi, \cdot), U, 0)=1.
$$
Note that condition (\ref{40}) holds true, for instance, if
$$
\langle \phi(t, x_0(r), y(t)), \, \dot x_0(r)\rangle
>0, \quad t,r\in [0,2\pi], \quad \|y\|_{C_T}\le \dfrac{1}{a},
$$
hence assumption (A$_3$) is also satisfied.

\vspace{2mm}
In conclusion, Theorem 1 applies to state the
existence, for $\e>0$ sufficiently small,  of a $2\pi$-periodic
solution $(x_\e, y_\e)$ to system (\ref{31}) such that
$$
\O(0,t,x_\e(t))\in U, \quad t\in [0,2\pi], \quad \|y_\e\|_{C_T}\le \dfrac{1}{a}.
$$
\qed
\begin{remark} We would like to point out that our result, as the proposed
example shows,
is not a consequence of known results concerning
the existence of periodic solutions around an equilibrium point with non zero
topological degree. In fact, if, for fixed $\xi\in U$, the function
$g(t,\xi)$ is not constant with respect to the time $t\in [0,2\pi]$,
then the differential equation (\ref{34}) does not have equilibrium points in
$U$.
\end{remark}

We end the paper by showing how from our result we can derive a classical result
of the theory of ordinary differential equations, cf. (\cite{3}, Theorem 3.1, p.
362).
Indeed, if we put in system (\ref{1}), $k=2, \, \psi_2(t,x)=0, A=0, y=0$, where
$x=(x_1,x_2)$, and $\psi_1(t,x)=Bx=(-x_2,x_1), \phi(t,x,0)=(0, g(t, -x_1, x_2)$,
where
$g$ is $2\pi$-periodic with respect to $t$. Then we have

$$
\eta(2\pi,s,\xi(a,\theta))-\eta(0,s,\xi(a,\theta))=
\left( \begin{array}{cc}
         \cos\theta & \sin\theta \\
         -\sin\theta & \cos\theta \end{array} \right) H(a,\theta),
$$
where $\xi(a,\theta)=(-a\cos\theta,a\sin\theta)$ and
$$
H(a,\theta)=\int_0^{2\pi}\left(\begin{array}{c}
(\sin\tau)f(\tau+\theta,a\cos\tau,-a\sin\tau) \\
(\cos\tau)f(\tau+\theta,a\cos\tau,-a\sin\tau)
\end{array} \right) d\tau.
$$
In \cite{3} it is assumed the following condition
$$
\mbox{det} | H'(a_0,{\theta}_0)|\not=0,
$$
for some $a_0\not=0$. It is now easy to see that this condition ensures that
$$
|\mbox{deg}(\eta(2\pi,s,\xi(\cdot))-\eta(0,s,\xi(\cdot)), V, 0)|=1.
$$
where $V$ is a sufficiently small open set containing $(a_0,\theta_0)$.

\end{document}